\documentclass[12pt,twoside,leqno]{article}
\usepackage{amsmath}
\usepackage{amssymb}
\usepackage{amsxtra}
\usepackage{amscd}
\usepackage{amsthm}
\usepackage[mathscr]{eucal}

\theoremstyle{plain}
\newtheorem{thm}[subsection]{Theorem}

\newtheorem{lem}[subsection]{Lemma}

\theoremstyle{definition}

\newtheorem{rem}[subsection]{Remark}
\newtheorem{para}[subsection]{}

\newenvironment{pf}{\proof[\proofname]}{\endproof}

\usepackage{fancyhdr} 
\usepackage{amssymb}
\usepackage{graphicx}

\begin{document}

\title{Spaces of ratios of schemes and \'etale cohomology}
\author
{Kazuya Kato, Chikara Nakayama, Sampei Usui}
\date{}

\maketitle

\renewcommand{\mathbb}{\bold}

\newcommand\Cal{\mathcal}
\newcommand\define{\newcommand}

\define\gp{\mathrm{gp}}%
\define\fs{\mathrm{fs}}%
\define\an{\mathrm{an}}%
\define\mult{\mathrm{mult}}%
\define\Ker{\mathrm{Ker}\,}%
\define\Coker{\mathrm{Coker}\,}%
\define\Hom{\mathrm{Hom}\,}%
\define\Ext{\mathrm{Ext}\,}%
\define\rank{\mathrm{rank}\,}%
\define\gr{\mathrm{gr}}%
\define\cHom{\Cal{Hom}}
\define\cExt{\Cal Ext\,}%

\define\cA{\Cal A}
\define\cB{\Cal B}
\define\cC{\Cal C}
\define\cD{\Cal D}
\define\cE{\Cal E}
\define\cF{\Cal F}
\define\cG{\Cal G}
\define\cH{\Cal H}
\define\cM{\Cal M}
\define\cO{\Cal O}
\define\cS{\Cal S}
\define\cT{\Cal T}
\define\cN{\Cal N}
\define\cR{\Cal R}
\define\cV{\Cal V}
\define\fF{\frak F}
\define\Dc{\check{D}}
\define\Ec{\check{E}}

\newcommand{\N}{{\mathbb{N}}}
\newcommand{\Q}{{\mathbb{Q}}}
\newcommand{\Z}{{\mathbb{Z}}}
\newcommand{\R}{{\mathbb{R}}}
\newcommand{\C}{{\mathbb{C}}}
\newcommand{\bN}{{\mathbb{N}}}
\newcommand{\bQ}{{\mathbb{Q}}}
\newcommand{\bF}{{\mathbb{F}}}
\newcommand{\bZ}{{\mathbb{Z}}}
\newcommand{\bP}{{\mathbb{P}}}
\newcommand{\bR}{{\mathbb{R}}}
\newcommand{\bC}{{\mathbb{C}}}
\newcommand{\bbQ}{{\bar \mathbb{Q}}}
\newcommand{\ol}[1]{\overline{#1}}
\newcommand{\too}{\longrightarrow}
\newcommand{\respect}{\rightsquigarrow}
\newcommand{\compatible}{\leftrightsquigarrow}
\newcommand{\upc}[1]{\overset {\lower 0.3ex \hbox{${\;}_{\circ}$}}{#1}}
\newcommand{\Gmlog}{\bG_{m, \log}}
\newcommand{\Gm}{\bG_m}
\newcommand{\ep}{\varepsilon}
\newcommand{\Spec}{\operatorname{Spec}}
\newcommand{\val}{{\mathrm{val}}} 
\newcommand{\n}{\operatorname{naive}}
\newcommand{\bs}{\operatorname{\backslash}}
\newcommand{\Gal}{\operatorname{{Gal}}}
\newcommand{\gal}{{\rm {Gal}}({\bar \Q}/{\Q})}
\newcommand{\galp}{{\rm {Gal}}({\bar \Q}_p/{\Q}_p)}
\newcommand{\gall}{{\rm{Gal}}({\bar \Q}_\ell/\Q_\ell)}
\newcommand{\wep}{W({\bar \Q}_p/\Q_p)}
\newcommand{\wel}{W({\bar \Q}_\ell/\Q_\ell)}
\newcommand{\Ad}{{\rm{Ad}}}
\newcommand{\BS}{{\rm {BS}}}
\newcommand{\even}{\operatorname{even}}
\newcommand{\End}{{\rm {End}}}
\newcommand{\odd}{\operatorname{odd}}
\newcommand{\GL}{\operatorname{GL}}
\newcommand{\np}{\text{non-$p$}}
\newcommand{\g}{{\gamma}}
\newcommand{\G}{{\Gamma}}
\newcommand{\Lam}{{\Lambda}}
\newcommand{\La}{{\Lambda}}
\newcommand{\lam}{{\lambda}}
\newcommand{\la}{{\lambda}}
\newcommand{\uL}{{{\hat {L}}^{\rm {ur}}}}
\newcommand{\uQp}{{{\hat \Q}_p}^{\text{ur}}}
\newcommand{\sel}{\operatorname{Sel}}
\newcommand{\dt}{{\rm{Det}}}
\newcommand{\Sig}{\Sigma}
\newcommand{\fil}{{\rm{fil}}}
\newcommand{\SL}{{\rm{SL}}}
\newcommand{\spl}{{\rm{spl}}}
\newcommand{\st}{{\rm{st}}}
\newcommand{\Isom}{{\rm {Isom}}}
\newcommand{\Mor}{{\rm {Mor}}}
\newcommand{\bg}{\bar{g}}
\newcommand{\id}{{\rm {id}}}
\newcommand{\cone}{{\rm {cone}}}
\newcommand{\al}{a}
\newcommand{\ChL}{{\cal{C}}(\La)}
\newcommand{\Image}{{\rm {Image}}}
\newcommand{\toric}{{\operatorname{toric}}}
\newcommand{\torus}{{\operatorname{torus}}}
\newcommand{\Aut}{{\rm {Aut}}}
\newcommand{\Qp}{{\mathbb{Q}}_p}
\newcommand{\barQp}{{\mathbb{Q}}_p}
\newcommand{\Qpur}{{\mathbb{Q}}_p^{\rm {ur}}}
\newcommand{\Zp}{{\mathbb{Z}}_p}
\newcommand{\Zl}{{\mathbb{Z}}_l}
\newcommand{\Ql}{{\mathbb{Q}}_l}
\newcommand{\Qlur}{{\mathbb{Q}}_l^{\rm {ur}}}
\newcommand{\F}{{\mathbb{F}}}
\newcommand{\eps}{{\epsilon}}
\newcommand{\epsLa}{{\epsilon}_{\La}}
\newcommand{\epsLaVxi}{{\epsilon}_{\La}(V, \xi)}
\newcommand{\epsOLaVxi}{{\epsilon}_{0,\La}(V, \xi)}
\newcommand{\Qplin}{{\mathbb{Q}}_p(\mu_{l^{\infty}})}
\newcommand{\otimesQplin}{\otimes_{\Qp}{\mathbb{Q}}_p(\mu_{l^{\infty}})}
\newcommand{\galFl}{{\rm{Gal}}({\bar {\Bbb F}}_\ell/{\Bbb F}_\ell)}
\newcommand{\gallur}{{\rm{Gal}}({\bar \Q}_\ell/\Q_\ell^{\rm {ur}})}
\newcommand{\galFF}{{\rm {Gal}}(F_{\infty}/F)}
\newcommand{\galFv}{{\rm {Gal}}(\bar{F}_v/F_v)}
\newcommand{\galF}{{\rm {Gal}}(\bar{F}/F)}
\newcommand{\epsVxi}{{\epsilon}(V, \xi)}
\newcommand{\epsOVxi}{{\epsilon}_0(V, \xi)}
\newcommand{\plim}{\lim_
{\scriptstyle 
\longleftarrow \atop \scriptstyle n}}
\newcommand{\sig}{{\sigma}}
\newcommand{\ga}{{\gamma}}
\newcommand{\del}{{\delta}}
\newcommand{\Vss}{V^{\rm {ss}}}
\newcommand{\Bst}{B_{\rm {st}}}
\newcommand{\Dpst}{D_{\rm {pst}}}
\newcommand{\Dcrys}{D_{\rm {crys}}}
\newcommand{\DdR}{D_{\rm {dR}}}
\newcommand{\Fin}{F_{\infty}}
\newcommand{\Kla}{K_{\lambda}}
\newcommand{\Ola}{O_{\lambda}}
\newcommand{\Mla}{M_{\lambda}}
\newcommand{\Det}{{\rm{Det}}}
\newcommand{\Sym}{{\rm{Sym}}}
\newcommand{\LaSa}{{\La_{S^*}}}
\newcommand{\cX}{{\cal {X}}}
\newcommand{\MHG}{{\frak {M}}_H(G)}
\newcommand{\tauMla}{\tau(M_{\lambda})}
\newcommand{\Fvur}{{F_v^{\rm {ur}}}}
\newcommand{\Lie}{{\rm {Lie}}}
\newcommand{\cL}{{\cal {L}}}
\newcommand{\cW}{{\cal {W}}}
\newcommand{\fq}{{\frak {q}}}
\newcommand{\cont}{{\rm {cont}}}
\newcommand{\SC}{{SC}}
\newcommand{\Om}{{\Omega}}
\newcommand{\dR}{{\rm {dR}}}
\newcommand{\crys}{{\rm {crys}}}
\newcommand{\hatSig}{{\hat{\Sigma}}}
\newcommand{\rdet}{{{\rm {det}}}}
\newcommand{\ord}{{{\rm {ord}}}}
\newcommand{\BdR}{{B_{\rm {dR}}}}
\newcommand{\BdRO}{{B^0_{\rm {dR}}}}
\newcommand{\Bcrys}{{B_{\rm {crys}}}}
\newcommand{\Qw}{{\mathbb{Q}}_w}
\newcommand{\barkappa}{{\bar{\kappa}}}
\newcommand{\cP}{{\Cal {P}}}
\newcommand{\cZ}{{\Cal {Z}}}
\newcommand{\oppLa}{{\Lambda^{\circ}}}
\newcommand{\add}{{{\rm {add}}}}
\newcommand{\red}{{{\rm {red}}}}
\newcommand{\Sp}{{{\rm {Sp}}}}
\newcommand{\GSp}{{{\rm {GSp}}}}
\newcommand{\nspl}{{{\rm {nspl}}}}
\newcommand{\triv}{{{\rm {triv}}}}
\newcommand{\Rep}{{{\rm {Rep}}}}
\newcommand\twoheaddownarrow{\mathrel{\rotatebox[origin=c]{90}{$\twoheadleftarrow$}}}

\newcommand{\br}{{\bold r}}
\newcommand{\mild}{{{\rm {mild}}}}
\newcommand{\MHS}{{{\rm MHS}}}
\newcommand{\LMH}{{{\rm LMH}}}
\newcommand{\pre}{{{\rm pre}}}
\newcommand{\sat}{{{\rm sat}}}
\newcommand{\ste}{{{\rm ste}}}
\newcommand{\SO}{{{\rm SO}}}
\newcommand{\et}{\mathrm{\acute e t}}
\newcommand{\loget}{\mathrm{log\acute e t}}
\newcommand{\ket}{\mathrm{k\acute e t}}
\newcommand{\pe}{\frak p}

\begin{abstract} 
We introduce the arithmetic topologies on the spaces of ratios for log schemes and prove the basic properties of their cohomologies. 
\end{abstract}

\renewcommand{\thefootnote}{\fnsymbol{footnote}}
\footnote[0]{MSC2020: Primary 14A21; Secondary 14A20, 14F20} 

\footnote[0]{Keywords: log geometry, fs log scheme, \'etale cohomology, ratio space}

\section*{Contents}

\noindent \S\ref{s:intro}. Introduction 

\noindent\S\ref{s:ZE}. Ratio spaces of log schemes, \'etale topology

\noindent\S\ref{s:prod}. Product of an arithmetic ratio space and a usual topological space

\noindent \S\ref{s:coh}. Invariance of Cohomologies

\noindent\S\ref{s:Betti}. GAGA and the change of base field

\noindent References

\section{Introduction} 
\label{s:intro}

Ratio spaces are introduced in the Hodge context 
to study the degeneration of mixed Hodge structures and
to give a new connection between the toroidal compactifications and the Borel--Serre compactification (\cite{KNU4}, \cite{KNU5}, \cite{KNU6}).

In this paper, we extend the theory of ratio spaces from the analytic context 
to the arithmetic context. We extend the \'etale cohomology theory of schemes to the \'etale cohomology theory of the ratio spaces. 
We will apply this to the theory of Shimura varieties (\cite{FKKN}). 

In Section \ref{s:ZE}, we introduce the definitions of
the ratio space of an fs log scheme, and its Zariski site, \'etale site, log \'etale site. 
After a preparation in Section \ref{s:prod},
we prove a basic invariant theorem in Section \ref{s:coh}. This is an algebraic analogue of the analytic result in \cite{KNUcon} that the cohomology of an fs log analytic space is the same as the cohomology of its ratio space. 
In Section \ref{s:Betti}, we prove a GAGA type theorem and a kind of base change theorem for subspaces of ratio spaces. 

\bigskip

\noindent Convention.
  For a sheaf $F$ on a site $Y$ and a morphism $X \to Y$ of sites, the pullback of $F$ on $X$ is denoted by $F_X$ or simply by $F$. 

  For an fs log scheme $S$, we denote by $S_{\triv}$ the open subscheme of $S$ consisting of all points at which the log structure is trivial.
  
\bigskip

  K.\ Kato was 
partially supported by NFS grants DMS 1601861 and DMS 2001182.
C.\ Nakayama was 
partially supported by JSPS Grants-in-Aid for Scientific Research (C) 21K03199 and (C) 26K06723.
S.\ Usui was 
partially supported by JSPS Grants-in-Aid for Scientific Research (C) 17K05200.

\section{Ratio spaces of log schemes, \'etale topology}\label{s:ZE}
Here we consider the spaces of ratios for log schemes, the Zariski topology, 
the \'etale topology, and the log \'etale topology of them, and state the main theorem. 

\begin{para}\label{ZE1}
  Let $S$ be an fs log scheme and assume that the log structure of $S$ comes from the Zariski site of $S$. 
  Then we can define the set $S_{[:]}$ exactly in the same way as in \cite{KNU4} 4.2 as follows.

  Let $S_{[:]}$ be the set of all pairs $(s,r)$, where $s\in S$ and $r \in
  R((M_S/\cO_S^\times)_{\bar s})$.
  Here $R((M_S/\cO_S^\times)_{\bar s})$ is the set of ratios, that is, the set of maps
$$r:(M_S/\cO_S^\times)_{\bar s} \times (M_S/\cO_S^\times)_{\bar s} \smallsetminus \{(1,1)\} \to [0,\infty]$$
satisfying the following conditions (i)--(iii).

\medskip

(i) $r(g,f)=r(f,g)^{-1}$.

(ii) $r(f,g)r(g,h)=r(f,h)$ if $\{r(f,g), r(g,h)\}\neq \{0,\infty\}$.

(iii)  $r(fg,h)=r(f,h)+r(g,h)$. 

\begin{para}
  Next we define the Zariski topology of $S_{[:]}$.
  Though the definition of the topology of $S_{[:]}$ given in the analytic situation \cite{KNU4} 4.2, which uses $\log|f|/\log|g|$, does not work for schemes, we can define the
  topology as follows. 

  We define the Zariski topology of $S_{[:]}$ as the weakest topology which satisfies the following condition:

 \medskip
 
 Let $U$ be an open set of $S$ and let $f,g$ be sections on $U$  of the log structure $M_S$. Let $C$ be the closed subset of $U$ consisting of all points $s$ such that at least one of the stalks $f_{\bar s}$ and $g_{\bar s}$ is not invertible in $M_{S,\bar s}$, let $C_{[:]}$ be the inverse image of $C$ in $S_{[:]}$ and let $r(f,g):C_{[:]}\to [0, \infty]$ be the map $(s,r)\mapsto r(f_{\bar s},g_{\bar s})$. Let $B$ be a closed subset of $[0, \infty]$ and let $B'$ be the inverse image of $B$ in $C_{[:]}$ under $r(f,g)$. Then  $U_{[:]} \smallsetminus B'$ is open in $S_{[:]}$.

\medskip

The projection $S_{[:]}\to S\; ;\; (s,r)\mapsto s$ is continuous.
\end{para}

  \end{para}
 
\begin{para}
  \label{R(P)}
  The ratio space $R(P)$ for a sharp fs monoid $P$ is a 
  compact Hausdorff and contractible 
  topological manifold with boundary.
  In fact, it is weakly contractible by \cite{KNUcon} Theorem 2.6.
  Since it is a compact Hausdorff topological manifold with boundary by
  \cite{KNU4} Lemma 4.2.7 and \cite{KNU4} Corollary 4.2.22, 
  it is an absolute neighborhood retract.
  Hence by Whitehead's theorem, it is contractible. 
\end{para} 

\begin{para}\label{ZE2} Example. Let $k$ be a field and let $S=\Spec(k[T_1,T_2])$ endowed with the log structure generated by $T_1$ and $T_2$. Let $p\in S$ be the point $(0,0)$. Then $S_{[:]}$ is obtained from $S$ by replacing $p$ with $[0, \infty]$, where
  $(p,r)$ corresponds to $r(T_2,T_1) \in [0,\infty]$. 
 
 An open set of $S_{[:]}$ for the Zariski topology is either an open set of $S\smallsetminus \{p\}$  or a set of the form $U_{[:]}\smallsetminus B$, where $U$ is an open set of $S$ containing $p$ and $B$ is a closed subset of $[0,\infty]$, satisfying the following conditions: If $U$ contains the prime ideal $(T_1)$, $B$ does not contain $\infty\in [0,\infty]$. 
  If $U$ contains the prime ideal $(T_2)$, $B$ does not contain $0\in [0, \infty]$. 
\end{para}
 
 \begin{para}\label{ZE3} Let $S$ be an fs log scheme. We define the \'etale site  $S_{[:],\et}$. 

An object of this site is a pair $(U, V)$, where $U$ is an object of the (classical, not log) \'etale site of $S$ such that the log structure of $U$ (which is the pullback log structure from $S$) comes from its Zariski site, and $V$ is an open set of $U_{[:]}$ for the Zariski topology (\ref{ZE2}). A morphism $(U,V)\to (U', V')$ is a morphism $U\to U'$ which induces $V\to V'$. A covering $((U_i, V_i) \to (U, V))_i$ is  a family of morphisms such that $\bigsqcup_i V_i \to V$ is surjective. 
\end{para}

 \begin{para}
Let $S$ be an fs log scheme. We define the kummer log \'etale (k\'et, for short) site  $S_{[:],\ket}$. 

An object of this site is a pair $(U, V)$, where $U$ is an object of the k\'et site of $S$ such that the log structure of $U$ comes from its Zariski site, and $V$ is an open set of $U_{[:]}$ for the Zariski topology (\ref{ZE2}). A morphism $(U,V)\to (U', V')$ is a morphism $U\to U'$ which induces $V\to V'$. A covering $((U_i, V_i) \to (U, V))_i$ is  a family of morphisms such that 
$\bigsqcup_i V_i \to V$ is surjective. 
\end{para}

\begin{para}
  Let $S$ be an fs log scheme and 
  assume that the log structure of $S$ comes from the Zariski site.
  Let $C$ be a locally closed subset of $S_{[:]}$.
  We define a sheaf of sets on $C$ with respect to the Zariski, \'etale, or k\'et topology as a sheaf on $U$ whose restriction to $U \smallsetminus C$ is a final object,
  where $U$ is an open set of $S_{[:]}$ containing $C$ as a closed subset.
  The sheaves on $C$ defines a topos, which is independent of the choice of $U$.
  Thus we have $H^m_*(C,\cdot)$ for $m\in \bZ$ and for $*=\emptyset, \et, \ket$.
\end{para}
 
The basic theorem of invariance of cohomologies of ratio spaces is as follows. 
This is proved in Section \ref{s:coh}.

\begin{thm}\label{ZE4} 
Let $S$ be an fs log scheme. 
Then the following holds. 

$(1)$  Assume that the log structure of $S$ comes from the Zariski site of $S$. 
Then we have isomorphisms  $H^m(S, F) \overset{\cong}\to H^m(S_{[:]}, F)$ for 
all $m$ and for all 
sheaves $F$ on $S$ of 
abelian groups.

$(2)$  Assume that locally on $S$, the dimension of $S$ is finite. 
  Then we have isomorphisms $H^m_{\et}(S, F) \overset{\cong}\to H^m_{\et}(S_{[:]}, F)$  for
  all $m$ and for all 
  sheaves $F$ on $S_{\et}$ of 
  abelian groups.

$(3)$ Assume that locally on $S$, the dimension of $S$ is finite.
  Assume also that locally on $S$, $S \smallsetminus S_{\triv}$ is of equi-characteristic.
  Then we have isomorphisms $H^m_{\ket}(S, F) \overset{\cong}\to H^m_{\ket}(S_{[:]}, F)$  for
  all $m$ and for all 
  sheaves $F$ on $S_{\ket}$ of 
  abelian groups.
\end{thm}

\section{Product of an arithmetic ratio space and a usual topological space}
\label{s:prod}
  In this section, we prove Theorem \ref{t:prod} below, which is used in the proof of Theorem \ref{ZE4} in the next section. 

  \begin{para}
    Let $X$ be an fs log scheme.
    Let $R$ be a topological space. 
    We define the \'etale (resp.\ k\'et) site $(X_{[:]} \times R)_{\et}$
  (resp.\ $(X_{[:]} \times R)_{\ket}$) as follows. 

  An object of this site is a pair $(U, V)$, where $U$ is an object of the \'etale (resp.\ k\'et) site of $X$ such that the log structure of $U$ comes from its Zariski site,
  and $V$ is an open set of the topological space $U_{[:]} \times R$.
  A morphism $(U,V)\to (U', V')$ is a morphism $U\to U'$ which induces $V\to V'$. A covering $((U_i, V_i) \to (U, V))_i$ is  a family of morphisms such that 
$\bigsqcup_i V_i \to V$ is surjective. 

  For example, let $S$ be a scheme endowed with the log structure defined by the homomorphism
  $\bN^2 \to \cO_S$ sending all non-trivial elements to zero. 
  Then $(S_{[:]})_{\et}$ 
  in Section \ref{s:ZE} is equivalent to
  $(S^{\mathrm {cl}} \times [0,\infty])_{\et}$, 
  where $S^{\mathrm {cl}}$ denotes the underlying scheme of $S$ with the trivial log structure. 
  \end{para}

  \begin{thm}\label{t:prod}
    Let $X$ be an fs log scheme and let $R$ be a
compact Hausdorff and contractible topological manifold with boundaries. 
Assume that locally on $X$, the dimension of $X$ is finite. 
    Let $F$ be a complex of abelian sheaves on the \'etale or the k\'et site of $X_{[:]}$ which is bounded below.
    Let $p: X_{[:]} \times R \to X_{[:]}$ be the projection of topoi.
    Then we have an isomorphism $F\overset{\cong}\to Rp_*p^*F$.
\end{thm}

  We prove this theorem together with another theorem as follows, which is an analogue of the smooth base change theorem.

  \begin{thm}\label{t:sm} Let $X$, $R$ and $F$ be as in Theorem $\ref{t:prod}$ and let $f: X \to Y$ be a strict morphism of fs log schemes.
    Then for the diagram 
$$    \begin{CD}
      X_{[:]} \times R @> f_{[:]}\times {\mathrm {id}}>> Y_{[:]} \times R \\
      @V p VV @VV p V \\
      X_{[:]}  @> f_{[:]}>> Y_{[:]}
      \end{CD}$$
of topoi, we have an isomorphism $p^*Rf_{[:]*}F \overset{\cong}\to R(f_{[:]}\times {\mathrm{id}})_*p^*F.$
  \end{thm}

 \begin{lem}\label{l:prod}
 Assume that Theorem $\ref{t:prod}$ is valid for $(X, F)$. Then Theorem $\ref{t:sm}$ is valid for this $F$ and for every $f:X\to Y$.
 \end{lem}
 
 \begin{pf}
   Since for any point of $R$ has a fundamental system of neighborhoods consisting of 
  compact Hausdorff and contractible topological manifolds with boundaries, this is seen by localizing $Y_{[:]}\times R$ by localizing both $Y_{[:]}$ and $R$. 
 \end{pf}

 \begin{para} We prove Theorem \ref{t:prod}.
   We may assume that $d=\dim X$ is finite.
   First we prove the case $d=0$.
    We may assume that $X$ is reduced (cf.\ \cite{I}). 
   Further, we may assume that $X$ is the spectrum of a (log) separably closed field.
   Then $X_{[:]}$ is regarded as the ratio space of the stalk of $M_X/\cO^{\times}_X$, which is a topological space.
  Hence the conclusion follows from the usual proper base change theorem in the category of topological spaces. 

   Assume now that the theorem is valid for $X$ with $\dim X<d$.
   We may assume that $X$ is local.
   Let $i:Z \to X$ be the strict closed immersion from the closed point and $j:U\to X$ the strict open immersion from the complement.
   By the induction hypothesis, the theorem is valid for $(U, j^*_{[:]}F)$.
   Then by applying Lemma \ref{l:prod} replacing $X \to Y$ there by $U \to X$, we can see that the theorem is valid for $(X, Rj_{[:]*}j^*_{[:]}F)$.
Hence the theorem for $(X,F)$ is equivalent to that for $(Z,Ri^!_{[:]}F)$, which is valid by the case $d=0$.
\end{para}

\begin{para} Theorem \ref{t:sm} follows from Theorem \ref{t:prod} and Lemma \ref{l:prod}. 
\end{para}

\section{Invariance of Cohomologies}\label{s:coh}
We prove Theorem \ref{ZE4}.

\begin{para}
  Here we give the proof of Theorem \ref{ZE4} (1) first and explain how to modify the arguments for Theorem \ref{ZE4} (2) and (3). 

  That is, we first prove the following. 
\end{para}

\begin{thm}
\label{t:VB}
  Let $S$ be an fs log scheme.
  Assume that the log structure of $S$ comes from the Zariski site of $S$. 
  Let $F$ be a complex of abelian sheaves on the Zariski site 
of $S$ which is bounded below.
  Let $\pi:S_{[:]}\to S$ be the canonical morphism. Then 
  we have
$$F \overset{\cong}\to R\pi_*\pi^*F.$$
\end{thm}

\begin{para}
  We prove this theorem.
  First, let $U=S_{\triv}$.
  Then we are reduced to the case $U=\emptyset$ as follows.

Let $j:U\to S$ and $j_{[:]}:U\to S_{[:]}$ be inclusion maps. 
Let $i:Z:=S \smallsetminus U \to S$ be the strict closed immersion from the complement.
\end{para}

\begin{lem}
\label{l:VB}
  $\pi^*Rj_*j^*F= Rj_{[:]*}j^*F$. 
\end{lem}
  
\begin{pf}
  Let $s' \in S_{[:]}$ and let $s \in S$ be its image. 
  By the definition of the topologies of $S_{[:]}$, the set
  $\{V' \times_{S_{[:]}}U\,|\,V'$ is an open 
  neighborhood of $s'$ in $S_{[:]}\}$ coincides with the set 
  $\{V \times_{S}U\,|\,V$ is an open 
  neighborhood of $s$ in $S\}$.
  This implies the desired equality.
\end{pf}

\begin{para}By this lemma, we have $R\pi_*\pi^*(Rj_*j^*F)= R\pi_*Rj_{[:]*}j^*F= Rj_*j^*F$.
That is, Theorem \ref{t:VB} is valid for $(S, Rj_*j^*F)$.
Hence Theorem \ref{t:VB} for $(S,F)$ is equivalent to that for $(Z,Ri^!F)$.
Thus we may and will assume that $U$ is empty. 
\end{para}

\begin{para}
  Before proceeding the proof, we make some preparations. 
  First we introduce the following notation. 
  For an fs log scheme $X$ charted by an fs monoid $P$ and for an element $\pe \in \Spec(P)$, define $V(\pe) \subset X$ as the inverse image of the closure of $\{\frak p\}$ in $\Spec(P)$ with respect to the induced map $X \to \Spec(P)$, and define a closed subset $C(\pe) \subset X_{[:]}$
as $\{(x,r)\,|\, x \in V(\pe), r \in R((M_X/\cO^{\times}_X)_{\bar x}), r(f,g)=\infty$ if $f
\in \frak p(M_X/\cO^{\times}_X)_{\bar x}$ and $g \not\in \frak p(M_X/\cO^{\times}_X)_{\bar x}\}$. 
(Cf.\ \cite{KNU6} Section 1.5.)

\end{para}

\begin{para}\label{R}
  Let $X$ be an fs log scheme charted by an fs monoid $P$.  
  Let $\pe_1, \ldots, \pe_r \in \Spec P$. 
  Then the following holds.

\medskip

(1) $\bigcap_{i=1}^r C(\frak p_i)$ is non-empty if and only if the set $\{\frak p_1, \dots, \frak p_r\}$ is totally ordered for the inclusion relation and, when we denote by $\frak p$ the maximum element among them, $V(\frak p)$ is non-empty.

(2) If $\emptyset \not= \frak p_1\subsetneq \dots \subsetneq \frak p_r$, $\bigcap_{i=1}^r C(\frak p_i)$ is homeomorphic 
to $V(\frak p_r)_{[:]}\times R$, where $V(\pe_r)$ is endowed with the log structure by $P\smallsetminus \pe_r$ and
$$R={\textstyle \prod}_{i=1}^{r}R(
Q_i/Q_i^{\times}),$$
where $Q_i=(P\smallsetminus \frak p_{i-1})_{(\frak p_i\smallsetminus \frak p_{i-1})}$ $(1 \le i \le r)$ with $\frak p_0=\emptyset$.
\end{para}

\begin{para}\label{Vietoris--Begle}
  In \ref{R} (2),
if $V(\pe_r)$ satisfies the conclusion of Theorem \ref{t:VB}, then 
$f: \bigcap_{i=1}^r \; C(\frak p_i) \to V(\frak p_r)$ 
satisfies $F\overset{\cong}\to Rf_*f^*F$ for any complex $F$ on $V(\pe_r)$ which is bounded below.
  In fact, let $V=V(\frak p_r)$.
  We want to prove that the composition $V_{[:]} \times R \to V_{[:]} \to V$ sends the pullback of $F$ to $F$.
  Since the map $V_{[:]}\times R\to V_{[:]}$ sends the pullback of $F$ to the pullback of $F$ by the proper base change theorem together with \ref{R(P)}, 
  the statement 
  is reduced to the conclusion of Theorem \ref{t:VB} for $V$. 
\end{para}

\begin{para}
  \label{pf:VB}
  We return to the proof of Theorem \ref{t:VB}. 
  Recall that $S_{\triv}$ is now empty.
  We may assume that there is an fs chart $Q\to M_S$.
  We use the induction on the rank of $Q^{\gp}/Q^\times$.
  If the rank is 0 or 1, we have $S_{[:]}=S$ and the statement is trivial. 
  We may assume that the underlying scheme of $S$ is 
  local.
  Let $s$ be its closed point and let $P\to M_S$ be an fs chart which induces $P/P^{\times}\cong (M_S/\cO_S^\times)_{\bar s}$.
  It is enough to show the equality at the stalk of $s$. 
  Since $S_{\triv}$ is empty, we have
  $S_{[:]}=\bigcup_{i=1}^r C(\frak p_i)$, where $\frak p_i$ runs over the set of non-empty prime ideals of $P$.
  We have an exact sequence
$$0\to F_{S_{[:]}} \to C^0 \to C^1 \to C^2\to \cdots.$$
Here $C^0= 
\bigoplus_i F_{C(\frak p_i)}$ and in general, $C^i = \bigoplus_J F_{C_J}$, where $J$ ranges over all subsets of $\{1,\dots, r\}$ of 
cardinality $i+1$ 
and $C_J$ is the intersection $\bigcap_{j\in J} C(\frak p_j)$.
For any $i$, since $P\smallsetminus \frak p_i$ is a proper face of $P$, 
the rank of $(P\smallsetminus \frak p_i)^{\gp}/(P \smallsetminus \pe_i)^{\times}$ is strictly smaller than the rank of $P^{\gp}/P^\times$.
  Hence the conclusion of Theorem \ref{t:VB} for $V(\pe_i)$ is valid by the induction hypothesis. 
This together with \ref{R} and \ref{Vietoris--Begle} implies that
$R\pi_*\pi^*F$ is isomorphic to the associated single complex of the double complex
$D^0\to D^1 \to D^2\to \cdots$.
Here $D^i=\bigoplus_J F_{V(\pe)}$, where $J$
ranges over all subsets of $\{1,\dots, r\}$ of 
cardinality $i+1$ such that $\{\frak p_j \;|\;j\in J\}$ is totally ordered and
$\pe$ is the maximum element of each $J$. 
Since there is the maximum element $P\smallsetminus P^{\times}$ among $\{\frak p_i\}$, the stalk at $s$ of this complex is
quasi-isomorphic to $F_s$, which proves Theorem \ref{t:VB}, and hence, Theorem \ref{ZE4} (1).
\end{para}

\begin{para}
\label{pf:et}
  We explain the proof for Theorem \ref{ZE4} (2) by describing how to modify the above proof. 
  It is enough to show the \'etale analogue of Theorem \ref{t:VB} under the assumption in Theorem \ref{ZE4} (2).
  The reduction to the case where $S_{\triv}$ is empty is the same including the \'etale analogue of Lemma \ref{l:VB}.
  In \ref{R} (2), 
$\bigcap_{i=1}^r C(\frak p_i)$ is still isomorphic 
to $V(\frak p_r)_{[:]}\times R$ in the sense that their \'etale sites coincide. 
In \ref{Vietoris--Begle}, $V_{[:]}\times R\to V_{[:]}$ sends the pullback of $F$ to $F$
by Theorem \ref{t:prod} instead of the proper base change theorem. 
In \ref{pf:VB}, we may assume that $S$ is strictly local instead of only being local.
Then the same argument works.
\end{para}

\begin{para}
\label{pf:ket}
  We explain the proof for Theorem \ref{ZE4} (3) by describing how to modify the proof for the \'etale case in \ref{pf:et}.
  It is enough to show the k\'et analogue of Theorem \ref{t:VB} under the assumption in Theorem \ref{ZE4} (3).
  We may assume that $S\smallsetminus S_{\triv}$ is of equi-characteristic. 
  We use the k\'et analogue of Lemma \ref{l:VB}.
  In \ref{pf:VB}, we may assume that $S$ is log strictly local. 
  The remaining problem is that in \ref{R} (2), 
$\bigcap_{i=1}^r C(\frak p_i)$ is no more equivalent 
to $V(\frak p_r)_{[:]}\times R$ in the sense that their k\'et sites do not necessarily coincide. 
  For, in \ref{Vietoris--Begle}, we have to show $F \overset \cong \to Rf_*f^*F$, where $F$ is on the k\'et site on $V(\pe_r)$ whose log structure is by $P$, not by $Q:=P \smallsetminus \pe_r$.
  For the sake of distinction, we denote the last $V(\pe_r)$ with $P$ by $V(\pe_r)'$. 
  Then 
  what we want is that 

\smallskip  
  
\noindent $(*)$ \qquad  $V(\pe_r)'_{[:]}  \times R \to V(\pe_r)'$ sends the pullback of $F$ on $V(\pe_r)'_{\ket}$ to $F$,

\smallskip  

\noindent  where the log structure of $V(\pe_r)'_{[:]}$ is by $P$ but the ratio in $V(\pe_r)'_{[:]}$ is with respect to $Q$.
Recall that the induction hypothesis that Theorem \ref{t:VB} holds for $V(\pe_r)$ (with the log by $Q$) is that $V(\pe_r)_{[:]} \to V(\pe_r)$ sends the pullback of $G$ to $G$ for any $G$ on $V(\pe_r)_{\ket}$.
  By this together with Theorem \ref{t:prod}, we know 

\smallskip  
  
\noindent $(**)$ \qquad  $V(\pe_r)_{[:]}  \times R \to V(\pe_r)$ sends the pullback of $G$ on $V(\pe_r)_{\ket}$ to $G$.

\smallskip  

  We reduce $(*)$ to $(**)$.   
  Consider the projection $P \to P/Q$.
  By assumption, $V(\pe_r)$ is equi-characteristic, say $p$.
  Let $\widetilde P$ be the pullback of $(P/Q)\otimes_{\bN}\bigcup_{p\nmid n}\frac 1{\,n\,} \bN$ by $P\otimes \bQ \to (P/Q)\otimes \bQ$.
Then the k\'et site of 
$V(\pe_r)' \otimes_{\bZ[P]} \bZ[\widetilde P]$ is equivalent to that of $V(\pe_r)$ (with $Q$).
Hence by the k\'et localizing with respect to $V(\pe_r)' \otimes_{\bZ[P]} \bZ[\widetilde P] \to V(\pe_r)'$,
$(*)$ reduces to $(**)$ as desired. 
\end{para}

\section{GAGA and the change of base field}
\label{s:Betti}
We prove the comparison theorem between the Betti and the \'etale cohomologies and the base change theorem with respect to the separably closed field extension.

\begin{para}
\label{C(pe)}
In this section, let $k$ be a field.
Let $S$ be an fs log scheme over $k$. 
Let $Y \subset S_{[:]}$ be a closed set. 
Assume that locally on $S$, the following condition is satisfied.

\smallskip

\noindent (i) \ \  There is an fs chart $P \to \cO_S$ such that
  $Y$ coincides with $\bigcup_{i=1}^n \bigcap _{j=1}^{n_i} C(\pe_{ij})$ for some $\pe_{ij} \in \Spec(P)$. 

\end{para}

\begin{thm}
  \label{t:GAGA}
  Assume $k=\bC$ and that $S$ is locally of finite type over $k$. 
  Then the following holds.
  
  $(1)$ $H^m_{\et}(Y,F)=H^m(Y^{\an},F^{\an})$ for any $m \ge0$ and any torsion sheaf $F$ of abelian groups on $S_{\et}$.

  $(2)$ $H^m_{\ket}(Y,F)=H^m(Y^{\an,\log},F^{\an})$ for any $m\ge0$ and any torsion sheaf $F$ of abelian groups on $S_{\ket}$.
\end{thm}

\begin{pf}
  We may assume that the condition \ref{C(pe)} (i) is satisfied by the Mayer--Vietoris exact sequence. 
  Let $C_i=\bigcap _{j=1}^{n_i} C(\pe_{ij})$ ($1 \le i \le m)$. 
  By using the resolution $ 0\to F_Y \to \bigoplus F_{C_i} \to \bigoplus F_{C_i \cap C_j}\to  \cdots $, we reduce to the case $m=1$.
  By \ref{R} (1), we may assume that $\{\pe_{1j}\}_j$ is totally ordered.
  Let $\pe$ be the maximum element among $\{\pe_{1j}\}_j$.

  To prove (1), by the \'etale analogue of \ref{R} (2) (cf.\ \ref{pf:et}) , we can replace $Y$ by $V(\frak p)_{[:]} \times R$, where $R$ is the product of some ratio spaces.
  By the arguments in \ref{pf:et}, 
  we can further replace $Y$ by $V(\frak p)$.
  Then the statement is classical.

To prove (2),  by the k\'et analogue of \ref{R} (2) (cf.\ \ref{pf:ket}), we can replace $Y$ by $V(\frak p)'_{[:]} \times R$, where $R$ is as above. 
Then we 
have $H^m_{\ket}(V(\pe)'_{[:]} \times R, F)=
H^m_{\ket}(V(\pe)', F)=H^m((V(\pe)')^{\an,\log}, F^{\an})
=H^m((V(\pe)'_{[:]})^{\an,\log} \times R, F^{\an})$, where the first equality is by $(*)$ in \ref{pf:ket} 
and the second one is by \cite{KN}.
   These equalities imply the statement after replacing $Y$ by $V(\frak p)'_{[:]} \times R$.
\end{pf}

\begin{thm}
  Let $Y$ be as in $\ref{C(pe)}$.
  Assume that $k$ is separably closed.
  Let $k'$ be a separably closed field extension of $k$.
  Then the following holds.
  
  $(1)$ $H^m_{\et}(Y,F)=H^m_{\et}(Y_{k'},F)$ for any $m \ge0$ and any torsion sheaf $F$ of abelian groups on $S_{\et}$.
  
  $(2)$ $H^m_{\ket}(Y,F)=H^m_{\ket}(Y_{k'},F)$ for any $m \ge0$ and any torsion sheaf $F$ of abelian groups on $S_{\ket}$.
\end{thm}

\begin{pf} By the same argument as in the proof of Theorem \ref{t:GAGA}, we can replace $Y$ by $V(\frak p)$ in (1) and by $V(\frak p)'$ in (2).
  Then (1) is classical.
  We prove (2) after replacing $Y$ by $V(\pe)'$.
  By \cite{I}, we may assume that both $k$ and $k'$ are algebraically closed. 
  Then the proof is completed by the smooth base change theorem in the log \'etale cohomology in \cite{N}.
\end{pf}

\begin{rem}
  The analogues of the results in Section \ref{s:coh} and Section \ref{s:Betti} hold for $H^0$ with a sheaf of sets as a coefficient.   It is reduced to
the results with abelian coefficients in these sections. 
  The authors expect that the results extend to the case for $H^1$ with non-commutative coefficients.
  \end{rem}

\noindent {\rm Kazuya KATO
\\
Department of Mathematics
\\
University of Chicago
\\
Chicago, Illinois, 60637, USA}
\\
{\tt kkato@math.uchicago.edu}

\bigskip

\noindent {\rm Chikara NAKAYAMA
\\
Department of Economics 
\\
Hitotsubashi University 
\\
2-1 Naka, Kunitachi, Tokyo 186-8601, Japan}
\\
{\tt c.nakayama@r.hit-u.ac.jp}

\bigskip

\noindent
{\rm Sampei USUI${}^{\sharp}$
\\
Graduate School of Science
\\
Osaka University}
\\
\footnote[0]{${\sharp}$ \ \ the deceased}

\begin{thebibliography}{99}

\bibitem{FKKN}
{\sc Fukaya, T., Kato, K., Kajiwara, T., Nakayama, C.},
\newblock
    {\em \'Etale cohomology groups of parabolic boundaries of Shimura varieties, I},
preprint. 

\bibitem{I}
  {\sc Vidal, I.},
  \newblock
      {\em Morphismes log \'etales et descente par hom\'eomorphismes universels},
      C.\ R.\ Acad.\ Sci.\ Paris S\'er.\ I.\ Math.\ {\bf 332} 
      (2001), 239--244.
      
\bibitem{KN}
{\sc Kato, K., Nakayama, C.},
\newblock
{\em Log Betti cohomology, log \'etale cohomology, and log de Rham cohomology of log schemes over $\C$},
Kodai Math.\ J.\ {\bf 22} (1999), 161--186. 

\bibitem{KNU4}
{\sc Kato, K., Nakayama, C., Usui, S.},
\newblock
{\em Classifying spaces of degenerating mixed Hodge structures,
IV$:$ The fundamental diagram}, 
Kyoto J.\ Math.\, {\bf 58} (2018), 289--426.

\bibitem{KNU6}
{\sc Kato, K., Nakayama, C., Usui, S.},
\newblock
{\em Classifying spaces of degenerating mixed Hodge structures, 
VI{\rm:} Log real analytic functions and log $C^{\infty}$ functions}, 
Yokohama Math.\ J.\ {\bf 70} (2024), 1--128.

\bibitem{KNUcon}
{\sc Kato, K., Nakayama, C., Usui, S.},
\newblock
{\em Toroidal compactifications and Borel--Serre compactifications}, 
Proc.\ Japan.\ Acad.\ {\bf 101}, Ser.\ A (2025), 13--18.

\bibitem{KNU5}
{\sc Kato, K., Nakayama, C., Usui, S.},
\newblock
{\em Classifying spaces of degenerating mixed Hodge structures,
V$:$ Extended period domains and algebraic groups}, 
preprint, submitted. arXiv:2107.03561.

\bibitem{N}
{\sc Nakayama, C.}, 
\newblock
{\em Logarithmic \'etale cohomology, II}, 
Advances in Math., {\bf 314} (2017), 663--725.
\end{thebibliography}
\end{document}